\documentclass[dvips,ps,preprint]{imsart}
\usepackage[psamsfonts]{amsfonts}
\usepackage{amsmath,amssymb,amsthm}
\usepackage[dvips]{graphicx}
\usepackage{eepic}
\usepackage{amsthm,amsmath}
\usepackage[square,comma]{natbib}
\RequirePackage[colorlinks]{hyperref}

\RequirePackage{hypernat}

\doi{10.1214/09-PS152}
\pubyear{2009}
\volume{6}
\issue{0}
\firstpage{34}
\lastpage{61}

\startlocaldefs

\newcommand{\lwo}{\colon \!\!}
\newcommand{\rwo}{\! \colon \!\!}
\newcommand{\wick}[1]{\lwo#1\rwo}

\newcommand{\Zbold} {{\mathbb Z}}
\newcommand{\Rbold} {{\mathbb R}}
\newcommand{\Pbold} {{\mathbb P}}
\newcommand{\Scal}   {\mathcal{S}}
\newcommand{\Lcal}   {\mathcal{L}}
\newcommand{\Gcal}   {\mathcal{G}}

\def\eqalign#1\enalign{
    \begin{align}#1\end{align}
    }
\newcommand{\nnb}   {\nonumber \\}
\newcommand{\Wcal}   {\mathcal{W}}

\renewcommand{\to} {\rightarrow}
\newcommand{\sumtwo}[2]{\sum_{ \mbox{ \scriptsize
    $\begin{array}{c}
                        {#1} \\ {#2}
                        \end{array} $ }
    }
}

\newcommand{\R}{\Rbold}
\newcommand{\Z}{\Zbold}

\newcommand{\C}{\mathbb{C}}
\newcommand{\1}{{\mathbb I}}

\newcommand{\cL}{{\cal L}}

\newcommand{\phib}{\bar\phi}
\newcommand{\psib}{\bar\psi}
\newcommand{\ci}{\underline{i}}

\newcommand{\Ex}{\mathbb{E}}

\newtheorem{theorem}{Theorem}[section]
\newtheorem{lemma}    [theorem] {Lemma}
\newtheorem{prop}    [theorem] {Proposition}
\newtheorem{rk} [theorem] {Remark}
\newtheorem{cor} [theorem] {Corollary}
\numberwithin{equation}{section}

\newcommand{\Gsrw}{G^{\,\rm srw}}
\newcommand{\Gsaw}{G^{\,\rm saw}}
\newcommand{\Gwsaw}{G^{\,\rm wsaw}}
\newcommand{\Gloop}{G^{\,\rm loop}}

\newcommand{\carL}{M}

\endlocaldefs

\begin{document}

\begin{frontmatter}

\title{ Functional integral representations for self-avoiding walk\thanksref{T1}}
\thankstext{T1}{This is an original survey paper.}
\runtitle{Functional integral representations for SAW}


\begin{aug}
 \author{\fnms{David C.} \snm{Brydges}\thanksref{t2}\ead[label=e1]{db5d@math.ubc.ca}}
\address{Department of Mathematics\\
University of British Columbia, Vancouver, BC V6T 1Z2, Canada\\ \printead{e1}}
\end{aug}

\begin{aug}
\author{\fnms{John Z.} \snm{Imbrie}\corref{}\ead[label=e2]{ji2k@virginia.edu}}
\address{Department of Mathematics\\
University of Virginia, Charlottesville, VA 22904-4137, U.S.A.\\ \printead{e2}}
\end{aug}
\medskip
\and

\begin{aug}
\author{\fnms{Gordon} \snm{Slade}\thanksref{t2}\ead[label=e3]{slade@math.ubc.ca}}
\address{Department of Mathematics\\
University of British Columbia, Vancouver, BC V6T 1Z2, Canada\\ \printead{e3}}

\thankstext{t2}{Supported in part by NSERC of Canada.}

\runauthor{D.C. Brydges et al.} 
\end{aug}

\begin{abstract}
We give a survey and unified treatment of functional integral
representations for both simple random walk and some self-avoiding
walk models, including models with strict self-avoidance, with weak
self-avoidance, and a model of walks and loops.  Our representation
for the strictly self-avoiding walk is new. The representations have
recently been used as the point of departure for rigorous
renormalization group analyses of self-avoiding walk models in
dimension 4. For the models without loops, the integral
representations involve fermions, and we also provide an
introduction to fermionic integrals. The fermionic integrals are in
terms of anticommuting Grassmann variables, which can be
conveniently interpreted as differential forms.
\end{abstract}



\received{\smonth{6} \syear{2009}}

\tableofcontents

\end{frontmatter}

\section{Introduction}
\label{sec:intro}

The use of random walk representations for
functional integrals in mathematical physics has a long history going back
to Symanzik \cite{Syma69}, who showed how such representations
can be used to study quantum field theories.
Representations of this type were exploited systematically in
\cite{ACF83,BFS83II,BFS82,Dynk83,FFS92}.
It is also possible to use such representations in reverse, namely to rewrite
a random walk problem in terms of an equivalent problem for a functional integral.

Our goal in this paper is to provide an introductory survey of
functional integral representations for some problems connected with
self-avoiding walks, with both strict and weak self-avoidance. In
particular, we derive a new representation for the strictly
self-avoiding walk. These representations have proved useful
recently in the analysis of various problems concerning
4-dimensional self-avoiding walks, by providing a setting in which
renormalization group methods can be applied.  This has allowed for
a proof of $|x|^{-2}$ decay of the critical Green function and
existence of a logarithmic correction to the end-to-end distance for
weakly self-avoiding walk on a 4-dimensional hierarchical lattice
\cite{BEI92,BI03c,BI03d}. It is also the basis for work in progress
on the critical Green function for weakly self-avoiding walk on
$\Z^4$ and a particular (spread-out) model of strictly self-avoiding
walk on $\Z^4$ \cite{BS10}.  In addition, the renormalization group
trajectory for a specific
model of weakly self-avoiding walk on $\Z^3$
(one with upper critical dimension $3+\epsilon$)
has
been constructed in \cite{MS08}, in this context. In this paper, we
explain and derive the representations, but we make no attempt to
analyze the representations here, leaving those details to
\cite{BEI92,BI03c,BI03d,BS10,MS08}.

The representations we will discuss can be divided into two classes:
purely bosonic, and mixed bosonic-fermionic.  The bosonic representations
will be the most familiar to probabilists, as they are in terms of ordinary
Gaussian integrals.  They represent simple random walks,
and also systems of self-avoiding
and mutually-avoiding walks and loops.

The mixed bosonic-fermionic representations
eliminate the loops, leaving only the
self-avoiding walk.  They involve Gaussian integrals with
anticommuting Grassmann variables.  A classic reference for Grassmann integrals
is the text by Berezin \cite{Bere66}, and there is a short introduction
in \cite[Appendix~B]{Salm99}.
Such integrals, although familiar in physics, are less so in probability theory.
It turns out, however, that these more exotic integrals share many features
in common with ordinary Gaussian integrals.
One of our goals is to provide a minimal introduction to these integrals,
for probabilists.

Representations for self-avoiding walks go back to an observation of
de Gennes \cite{Genn72}.
The $N$-vector model has a random walk representation
given by a self-avoiding walk in a background of mutually-avoiding
self-avoiding loops, with
every loop contributing a factor $N$.
This led de Gennes to consider the limit $N \to 0$,
in which closed loops no longer contribute, leading to a representation
for the self-avoiding walk model
as the $N=0$ limit of the $N$-vector model (see also \cite[Section~2.3]{MS93}).
Although this idea has been very useful in physics, it has been less productive
within mathematics, because $N$ is a natural number and so it is unclear how
to understand a limit $N \to 0$ in a rigorous manner.

On the other hand, the notion was developed in \cite{McKa80,PS80}
that while an $N$-component boson field $\phi$ contributes a factor
$N$ to each closed loop, an $N$-component fermion field $\psi$
contributes a complementary factor $-N$. The net effect is to
associate zero to each closed loop. We give a concrete demonstration
of this effect in Section~\ref{sec:expres} below. This provides a
way to realize de Gennes' idea, without any nonrigorous limit.

Moreover, it was
pointed out by Le Jan \cite{LeJa87,LeJa88}
that the anticommuting variables can be represented by differential forms:
the fermion field can be regarded as nothing more than the differential
of the boson field.
This observation was further developed in \cite{BJS03,BI03c}, and we
will follow the approach based on differential forms
in this paper.
In this approach,
the anticommuting nature of fermions is
represented by the anticommuting wedge product for differential forms.
Thus the world of Grassmann variables, initially mysterious, can be
replaced by differential forms, objects which are fundamental
in differential geometry in the way that random variables are fundamental
in probability.

We have attempted
to keep this paper self-contained.  In particular,
our discussion of differential forms for the representations involving
fermions is intended to be introductory.

The rest of the paper is organized as follows. In
Section~\ref{sec:wlm}, we derive integral representations for simple
random walk, and for a model of a self-avoiding walk and
self-avoiding loops all of which are mutually avoiding.  These are
purely bosonic representations, without anticommuting fermionic
variables. In Section~\ref{sec:sawm}, we define the self-avoiding
walk models (without loops).  Their representations are derived in
Section~\ref{sec:fermion}, using the fermionic integration
introduced in Section~\ref{sec:rwdf}. The mixed bosonic-fermionic
integrals are examples of \emph{supersymmetric} field theories.
Although an appreciation of this fact is not necessary to understand
the representations, in Section~\ref{sec:susy} we briefly discuss
this important connection.

\section{Bosonic representations}
\label{sec:wlm}

\subsection{Gaussian integrals}

By ``bosonic representations'' we mean representations for random walk
models in terms of ordinary Gaussian integrals.
For our purposes, these integrals are in terms of a two-component
field $(u_x,v_x)_{x\in \{1,\ldots,\carL\}}$,
which is most conveniently represented by the complex pair $(\phi_x,\phib_x)$,
where
\begin{equation}
    \phi_x = u_x + i v_x, \quad \phib_x = u_x - iv_x.
\end{equation}
The differentials $d\phi_x$, $d\phib_x$ are given by
\begin{equation}
    d\phi_x = du_x + idv_x, \quad d\phib_x = du_x - i dv_x,
\end{equation}
and their product $d\phib_x d\phi_x$ is given by
\begin{equation}
\label{e:ddphi}
    d\phib_x d\phi_x
    = 2 i du_x dv_x,
\end{equation}
where we adopt the convention that differentials are multiplied together
with the anticommutative wedge product; in particular $du_x du_x$
and $dv_x dv_x$ vanish and do not appear in the above product.
This anticommutative product will play a central role when we come to
fermions in Section~\ref{sec:rwdf}, but until then plays no role beyond
the formula \eqref{e:ddphi}.  We are using the letter ``$x$'' as index
for the field in anticipation of the fact that in our representations
the field will be indexed by the space in which our random walks take
steps.

We now briefly review
some elementary properties of Gaussian measures.
Let $C$ be an $\carL \times \carL$ complex matrix.
We assume that $C$ has positive Hermitian part,
i.e., $\sum_{x,y=1}^\carL \phi_x (C_{x,y}+ \bar{C}_{y,x})\phib_y >0$ for all
nonzero $\phi \in \C^\carL$.
Let $A=C^{-1}$.  We write $d\mu_C$ for the Gaussian measure
on $\R^{2\carL}$ with covariance $C$, namely
\begin{equation}
    d\mu_C(\phi,\phib)
    =
    \frac{1}{Z_C} e^{-\phi A\phib}
    d\phib_1 d\phi_1 \cdots
    d\phib_\carL d\phi_\carL
    ,
\end{equation}
where $\phi A \phib = \sum_{x,y=1}^\carL \phi_x A_{x,y}\phib_y$, and
where $Z_C$ is the normalization constant
\begin{equation}
    Z_C = \int_{\R^{2\carL}} e^{-\phi A\phib}
    d\phib_1 d\phi_1 \cdots
    d\phib_\carL d\phi_\carL.
\end{equation}
We will need the value of $Z_C$ given in the following lemma.

\begin{lemma}
\label{lem:ZC}
For $C$ with positive Hermitian part and inverse $A=C^{-1}$,
\begin{equation}
\label{e:Geval}
    Z_C
    =
    \int e^{-\phi A \bar\phi}
    d\phib_1 d\phi_1 \cdots
    d\phib_\carL d\phi_\carL
    =
    \frac{(2\pi i)^\carL}{\det A}.
\end{equation}
\end{lemma}

\begin{proof}
Consider first the case where $C$, and hence $A$, is Hermitian.
In this case, there is a unitary matrix $U$ and a diagonal
matrix $D$ such that $A=U^{-1}DU$.
Then $\phi A \phib = w D \bar{w}$, where $w = U\phi$, so
\begin{equation}
    \frac{1}{(2\pi i)^\carL} Z_C
    =
    \prod_{x=1}^\carL
    \left(
    \frac{1}{\pi}
    \int_{-\infty}^\infty e^{-d_x (u_x^2+v_x^2)}
    du_xdv_x  \right)
    = \prod_{x=1}^\carL\frac{1}{d_x} = \frac{1}{\det A}.
\end{equation}
For the general case, we write $A (z) =G+iz H$ with $G= \frac 12
(A+A^\dagger)$, $H=\frac{1}{2i}(A-A^\dagger)$ and $z=1$.  Since
$\phi (iH) \phib$ is imaginary, when $G$ is positive definite the
integral in \eqref{e:Geval} converges and defines an analytic
function of $z$ in a neighborhood of the real axis.  Furthermore,
for $z$ small and purely imaginary, $A (z)$ is Hermitian and
positive definite, and hence \eqref{e:Geval} holds in this case.
Since $(\det A(z))^{-1}$ is a meromorphic function of $z$,
\eqref{e:Geval} follows from the uniqueness of analytic extension.
\end{proof}

A basic tool is the integration by
parts formula given in the following lemma.
The derivative appearing in its statement
is defined by
\begin{equation}
\label{e:ddphix}
    \frac{\partial }{\partial \phi_x}
    =\frac 12 \left(\frac{\partial }{\partial u_x}
    - i \frac{\partial }{\partial v_x}\right).
\end{equation}
With $\partial /\partial \phib_x$ defined to be its conjugate,
this leads to the equations
\begin{equation}
\label{e:phiderivs}
    \frac{\partial \phi_y}{\partial \phi_x}
    =
    \frac{\partial \phib_y}{\partial \phib_x} = \delta_{x,y},
    \quad\quad
    \frac{\partial \phib_y}{\partial \phi_x}
    =
    \frac{\partial \phi_y}{\partial \phib_x} = 0.
\end{equation}

\begin{lemma}
\label{lem:ibp}
Let $C$ have positive Hermitian part. Then
\begin{equation}
\label{e:Gibp}
    \int_{\R^{2\carL}} \phib_a F \, d\mu_C(\phi,\phib)
    =
    \sum_{x\in \Lambda} C_{a,x} \int_{\R^{2\carL}}
    \frac{\partial F}{\partial \phi_x}
    \, d\mu_C(\phi,\phib),
\end{equation}
where $F$ is any $C^1$ function such that both sides are integrable.
\end{lemma}

\begin{proof}
Let $A=C^{-1}$.
We begin with the integral on the right-hand side,
and make the abbreviation
$d\phib d\phi = d\phib_1d\phi_1 \cdots d\phib_\carL d\phi_\carL$.
By \eqref{e:ddphix}, we can use
standard integration by parts to move the derivative from one factor
to the other, and with \eqref{e:phiderivs} this gives
\begin{equation}
    \int
    \frac{\partial F}{\partial \phi_x}
    e^{-\phi A \phib} d\phib d\phi
    =
    -\int
    \frac{\partial e^{-\phi A \phib}}{\partial \phi_x} F
    d\phib d\phi
    =
    \int \sum_{y} A_{x,y}\phib_y F e^{-\phi A \phib} d\phib d\phi.
\end{equation}
Now we multiply by $C_{a,x}$, sum over $x$, and use $C=A^{-1}$,
to complete the proof.
\end{proof}

The equations
\eqalign
\label{e:cov}
    \int_{\R^{2\carL}} \phi_a \phi_b \, d\mu_C(\phi,\phib) &=
    \int_{\R^{2\carL}} \phib_a \phib_b \, d\mu_C(\phi,\phib) =
    0,
    \nnb
    \int_{\R^{2\carL}} \phib_a \phi_b \, d\mu_C(\phi,\phib) &= C_{a,b}.
\enalign
are simple consequences of Lemma~\ref{lem:ibp}.  The last equality
is a special case of Wick's theorem, which provides a formula
for the calculation of arbitrary moments of the Gaussian measure.
We will only need the following special case of Wick's theorem, in
which a particular Gaussian expectation is evaluated as the permanent
of a submatrix of $C$.

\begin{lemma}
\label{lem:Wick}
Let $\{x_1,\ldots,x_k\}$ and $\{y_1,\ldots,y_k\}$ each be
sets of $k$ distinct points in $\Lambda$, and let $S_k$ denote the
set of permutations of $\{1,\ldots,k\}$.
Then
\begin{equation}
    \int_{\R^{2\carL}}
    \left(\prod_{l=1}^k \phib_{x_l}\phi_{y_l} \right)
    d\mu_C(\phi,\phib)
    =
    \sum_{\sigma \in S_k} \prod_{l=1}^k C_{x_l,\sigma(y_l)}.
\end{equation}
\end{lemma}

\begin{proof}
This follows by repeated use of integration by parts.
\end{proof}


\subsection{Simple random walk}

Our setting throughout the paper
is a fixed finite set $\Lambda=\{1,2,\ldots,\carL\}$
of cardinality $\carL \ge 1$.  Given points $a,b \in \Lambda$,
a \emph{walk} $\omega$
from $a$ to $b$ is a sequence of points
$x_0=a,x_1,x_2,\ldots, x_n=b$, for some $n \ge 0$.
We write $|\omega|$ for the length $n$ of $\omega$.
Sometimes it is useful to regard $\omega$ as consisting of the
directed edges $(x_{i-1},x_i)$, $1 \le i \le n$, rather than vertices.
Let $\Wcal_{a,b}$ denote the set of all walks from $a$ to $b$, of any length.

Let $J$ be a $\Lambda \times \Lambda$ complex matrix with zero diagonal
part (i.e., $J_{x,x}=0$ for all $x \in \Lambda$).
Let $D$ be a diagonal matrix with nonzero entries $D_{x,x}=d_x \in \C$.
We assume that $D-J$ is
\emph{diagonally dominant}; this means that
\begin{equation}
    \max_{x \in \Lambda}
    \sum_{y \in \Lambda} \left| \frac{J_{x,y}}{d_{x}}\right|
    < 1.
\end{equation}
Given $\omega \in \Wcal_{a,b}$, let
\begin{equation}
    J^w = \prod_{e\in\omega} J_{e}.
\end{equation}
Here we regard $\omega$ as a set
of labeled edges $e = (\omega(i-1),\omega(i))$
(the empty product is $1$ if $|\omega|=0$).
The simple random walk two-point function is defined by
\eqalign
\label{e:Gsrw}
    \Gsrw_{a,b}
    & =
    \sum_{\omega \in \Wcal_{a,b}} J^\omega
    \prod_{i=0}^{|\omega|} d_{\omega(i)}^{-1}.
\enalign
The assumption that $D-J$ is diagonally dominant ensures that
the sum in \eqref{e:Gsrw} converges absolutely.
The following theorem was proved in \cite{BFS82}.

\begin{theorem}
\label{thm:srw}
Suppose that $D-J$ is diagonally dominant.  Then
$C=(D-J)^{-1}$ exists and
$\Gsrw_{a,b} = (D-J)^{-1}_{a,b}$.
In addition, if  $D-J$ has positive Hermitian part then
\begin{equation}
    \Gsrw_{a,b}
    =
    (D-J)^{-1}_{a,b}
    = \int_{\R^{2\carL}} \phib_a \phi_b \, d\mu_C(\phi,\phib).
\end{equation}
\end{theorem}

\begin{proof}
The sum in \eqref{e:Gsrw} can be evaluated explicitly as
\begin{equation}
\label{e:DJDa}
    \Gsrw_{a,b}
    =
    \sum_{\omega \in \Wcal_{a,b}} J^\omega
    \prod_{i=0}^{|\omega|} d_{\omega(i)}^{-1}
    =
    \sum_{n=0}^\infty \left( D^{-1} (JD^{-1})^n ) \right)_{a,b}
    .
\end{equation}
It is easily verified that $D-J$ applied to the right-hand side gives
the identity, and hence
\begin{equation}
\label{e:DJD}
    \Gsrw_{a,b}
    =
    (D-J)^{-1}_{a,b}.
\end{equation}
When $D-J$ has positive Hermitian part, we may use \eqref{e:cov} to
complete the proof.
\end{proof}

Next, we suppose that $d_x>0$, $J_{x,y} \ge 0$, and
give two alternate representations for $\Gsrw_{a,b}$ in
terms of continuous-time Markov chains.  For the first, which appeared
in \cite{Dynk83}, we consider the continuous-time Markov chain
$X$ defined as follows.  The state space of $X$ is $\Lambda \cup
\{\partial\}$, where $\partial$ is an absorbing state called the
cemetery.  When $X$ arrives at state $x$ it waits for an
${\rm Exp}(d_x)$ holding time and then jumps to $y$ with probability
$\pi_{x,y}=d_x^{-1}J_{x,y}$ and jumps to the cemetery with probability
$\pi_{x,\partial} = 1- \sum_{y \in \Lambda}d_x^{-1}J_{x,y}$.
The holding times are
independent of each other and of the jumps.
Let $\zeta$ denote the time at which the process arrives in the
cemetery.  Note that
if $D-J$ is diagonally dominant then $\zeta < \infty$ with probability $1$,
and by right-continuity of the sample paths the last state visited by
$X$ before arriving in the cemetery is $X(\zeta^-)$.
For $x \in \Lambda$, let $L_x$ denote the
total (continuous) time spent by $X$ at $x$.
We denote the expectation for $X$, started from $a \in \Lambda$,
by $\Ex_a$.

\begin{theorem}
\label{thm:srw2}
Suppose that $D-J$ is diagonally dominant, with $d_x>0$, $J_{x,y}\ge 0$,
and let $\overline d_x = \sum_{y\in\Lambda}J_{x,y}$.
Let $V$ be a diagonal matrix with entries $V_{x,x}=v_x$, and
suppose that $0<\overline d_x < d_x + {\rm Re}\,v_x $ for all
$x \in \Lambda$.
Let $\Gsrw_{a,b}(v)$ denote the two-point function \eqref{e:Gsrw},
with matrix $D+V-J$ in place of $D-J$.  Then
\begin{equation}
\label{e:Dynkin}
    \Gsrw_{a,b}(v)
    =
    \frac{1}{d_b\pi_{b,\partial}}
    \Ex_a \left( e^{-\sum_{x\in \Lambda}v_x L_x}
    \1_{X(\zeta^-) = b} \right).
\end{equation}
\end{theorem}

\begin{proof}
The Markov chain $X$
is equivalent to a discrete-time Markov chain $Y$ which jumps
with the above transition probabilities, together with a sequence
$\sigma_0,\sigma_1,\ldots$ of exponential holding times.
Let $\eta$ denote the discrete random time after
which the process $Y$
jumps to $\partial$.
By partitioning on the events $\{\eta =n\}$, noting that $\eta$ is
almost surely finite, we see that the right-hand side of
\eqref{e:Dynkin} is equal to
\begin{equation}
\label{e:Dynkin1}
    \frac{1}{d_b}\sum_{n=0}^\infty \Ex_a \left(
    e^{-\sum_{i=0}^{n}v_{Y_i}\sigma_i}
    \1_{Y_n = b} \right).
\end{equation}
Given the sequence $Y_0,Y_1,\ldots,Y_n$, the $\sigma_i$ are independent
${\rm Exp}(d_{Y_i})$ random variables and hence
\begin{equation}
    \frac{1}{d_b}\Ex_a\left( e^{-\sum_{i=0}^{n}v_{Y_i}\sigma_i}
    \1_{Y_n = b} | Y_0,Y_1,\ldots \right)
    =
    \frac{1}{d_b+v_b}\prod_{i=0}^{n-1} \frac{d_{Y_i}}{d_{Y_i}+v_{Y_i}}
    .
\end{equation}
If we then take the expectation with respect to the Markov chain $Y$,
we find that \eqref{e:Dynkin1} is equal to
\begin{equation}
    \sum_{n=0}^\infty \sum_{\omega \in \Wcal_{a,b}: |\omega|=n}
    \pi^\omega
    \frac{1}{d_b+v_b}
    \prod_{i=0}^{n-1} \frac{d_{\omega(i)}}{d_{\omega(i)}+v_{\omega(i)}}
    =
    \sum_{\omega \in \Wcal_{a,b}} J^\omega
    \prod_{i=0}^{|\omega|} \frac{1}{d_{\omega(i)}+v_{\omega(i)}},
\end{equation}
which is the desired result.
\end{proof}

Next, we derive a third representation for $\Gsrw_{a,b}(v)$,
which is more general than Theorem~\ref{thm:srw2} as it does not
require diagonal dominance of $D-J$ (it does require
${\rm Re} \,v_x>0$ when $d_x=\overline d_x$).
This representation was obtained in \cite{BEI92} using
the Feynman--Kac formula,
but we give a different proof based on Theorem~\ref{thm:srw2}.
The representation involves a second continuous-time
Markov process, with
generator $\overline D - J$ where we set
$\overline d_x= \sum_{y \in \Lambda}J_{x,y}$
and assume $\overline d_x >0$ for each $x \in \Lambda$.
This process is like the one described above, but has no
cemetery site and  continues for all time.
Let $\overline\Ex_a$ denote the expectation for this process started at
$a\in \Lambda$.
Let
\begin{equation}
\label{e:LTdef}
    L_{x,T} = \int_0^T \1_{X(s)=x} ds.
\end{equation}
denote the time spent by $X$
at $x$ during the time interval $[0,T]$.

%
%

\begin{theorem}
\label{thm:srw3}
Suppose that $d_x>0$, $J_{x,y}\ge 0$,
and let $\overline d_x = \sum_{y\in\Lambda}J_{x,y}$.
Let $V$ be a diagonal matrix with entries $V_{x,x}=v_x$, and
suppose that $0<\overline d_x < d_x + {\rm Re}\,v_x $ for all
$x \in \Lambda$.
Then
\begin{equation}
\label{e:FKrep}
    \Gsrw_{a,b}(v)
    =
    \int_0^\infty \overline\Ex_a
    \left( e^{-\sum_{x\in \Lambda}(v_x+d_x - \overline d_x)L_{x,T}}
    \1_{X(T)=b} \right)
    dT.
\end{equation}
\end{theorem}

\begin{proof}
Let $\mu = \min_{x\in\Lambda} ({\rm Re}\,v_x + d_x - \overline d_x)$
and let
$0 < \epsilon < \mu$.  We write
\begin{equation}
    D+V-J = D^{(\epsilon)} + V^{(\epsilon)} - J
\end{equation}
with
\begin{equation}
    D^{(\epsilon)}_{x,x}
    =
    d_x^{(\epsilon)}
    =
    \overline d_x + \epsilon,
    \quad
    V^{(\epsilon)}_{x,x}
    = v_x^{(\epsilon)}
    = v_x + d_x- \overline d_x - \epsilon .
\end{equation}
Let $\Ex_{a}^{(\epsilon)}$ denote the expectation for
the Markov process defined in terms of
$D^{(\epsilon)} - J$.
Since $D^{(\epsilon)}  - J$ is diagonally dominant
and ${\rm Re}\, v_x^{(\epsilon)}\ge \mu-\epsilon$, by
Theorems~\ref{thm:srw} and \ref{thm:srw2} we have
\eqalign
    \Gsrw_{a,b} (v)
    &=
    (D+V-J)^{-1}_{a,b}
    = (D^{(\epsilon)} + V^{(\epsilon)} - J)^{-1}_{a,b}
    \nnb &
    =
    \frac{1}{\epsilon}
    \Ex^{(\epsilon)}_a \left(
    e^{-\sum_{x\in \Lambda}v_x^{(\epsilon)}L_{x}}
    \1_{X(\zeta^-) =b} \right)
    ,
\enalign
where the $\epsilon$ in the
denominator is equal to the product of $d_b^{(\epsilon)}$ and
$\pi^{(\epsilon)}_{b,\partial}=\epsilon/d_b^{(\epsilon)}$.

We
partition on the values of $\zeta$, the time of transition to
$\partial$. For $\delta >0$, let
\begin{equation}
    I (\delta )
    =
    \{j \delta :j = 0,1,2,\dotsc \}.
\end{equation}
Then
\begin{equation}
    \Gsrw_{a,b} (v)
    =
    \sum_{T \in I (\delta)}
    \frac{1}{\epsilon}
    \Ex^{(\epsilon)}_a \left(
    e^{-\sum_{x\in \Lambda}v_x^{(\epsilon)}L_{x}}
    \1_{Y_\eta =b} \1_{T <\zeta \le T+\delta}\right).
\end{equation}
The probability of the symmetric difference
\begin{equation}
    \{Y_\eta =b, T <\zeta \le T + \delta\}
    \Delta
    \{X(T)=b,X(T + \delta)=\partial \}
\end{equation}
is $O (\delta^{2})$ because this event requires two jumps in time
$\delta $.
Also, $L_{x,T} \le L_{x} \le L_{x,T}+\delta$ on the event $\{T
<\zeta \le T + \delta \}$, so
\begin{equation}
    \Gsrw_{a,b} (v)
    =
    \lim_{\delta \rightarrow 0}
    \sum_{T \in I (\delta )}
    \frac{1}{\epsilon}
    \Ex^{(\epsilon)}_a \left(
    e^{-\sum_{x\in \Lambda}v_x^{(\epsilon)} L_{x,T}}
    \1_{X(T)=b,X(T + \delta)=\partial }\right).
\end{equation}
By the Markov property and the fact that
\begin{equation}
    \overline{\Pbold} (X(T + \delta)=\partial |X(T)=b)
    =
    d_b^{(\epsilon)}\delta
    \pi^{(\epsilon)}_{b,\partial} + O (\delta^{2})
    =\epsilon\delta + O(\delta^2),
\end{equation}
we obtain
\begin{align}
    \Gsrw_{a,b} (v)
    &=
    \lim_{\delta \rightarrow 0}
    \sum_{T \in I (\delta )}
    \Ex^{(\epsilon)}_a \left(
    e^{-\sum_{x\in \Lambda}v_x^{(\epsilon)}L_{x,T}}
    \1_{X(T)=b}\right) \delta
    \nnb
    &=
    \int_{0}^{\infty }
    \Ex^{(\epsilon)}_a \left(
    e^{-\sum_{x\in \Lambda}v_x^{(\epsilon)}L_{x,T}}
    \1_{X(T)=b}\right) \,dT .
\end{align}
Now taking the limit $\epsilon \to 0$,
$\Ex^{(\epsilon)}_a$ converges to $\overline\Ex_a$ on
bounded functions of $\{X(t) : 0 \leq t \leq T\}$ since the transition
probabilities and the densities of the holding times $\sigma_i$
converge to their analogues in $\overline\Ex_a$.  Noting that
\begin{equation}
    \left|
    e^{-\sum_{x\in \Lambda}v_x^{(\epsilon)} L_{x,T}}
    \right|
    \le
    e^{-(\mu-\epsilon)T},
\end{equation}
we obtain \eqref{e:FKrep} by dominated convergence.
\end{proof}

The two representations for $\Gsrw_{a,b}$ in
Theorems~\ref{thm:srw2}--\ref{thm:srw3} show that the right-hand sides
of \eqref{e:Dynkin} and \eqref{e:FKrep} are equal.
 The following
 proposition generalizes this equality.

\begin{prop}
\label{prop:DBEI} Suppose that $D-J$ is diagonally dominant, with
$d_x>0$, $J_{x,y}\ge 0$.
Fix $0<\epsilon < \min_{x\in \Lambda} (d_x - \overline d_x)$.
Let $F: [0,\infty)^{\carL } \to
\C$ be a Borel function such that there is a constant $C$ for which
$|F (t)| \le C\exp( \epsilon \sum_{x}t_{x})$.
Let $L = (L_{x})_{x \in \Lambda }$ and similarly for $L_T$.
Then
\begin{equation}
\label{e:DBEI}
    \frac{1}{d_b\pi_{b,\partial}}
    \Ex_a \left(
    F (L)
    \1_{X(\zeta^-) = b} \right)
    =
   \int_{0}^{\infty }
    \overline\Ex_a \left(
    F (L_{T})
    e^{-\sum_{x\in \Lambda}(d_x - \overline d_x )L_{x,T}}
    \1_{X(T)=b}\right) \,dT .
\end{equation}
\end{prop}

\begin{proof}
Let $S$ be a Borel subset of $[0,\infty)^{\carL }$,
and let $\chi_S$ denote the characteristic function of $S$.
We define $\mu(S)$ and $\nu(S)$ by evaluating the
left- and right-hand sides of \eqref{e:DBEI} on $F=\chi_S$, respectively.
With these definitions, $\mu$ and $\nu$ are finite Borel measures.
Together, Theorems~\ref{thm:srw2}--\ref{thm:srw3} establish
\eqref{e:DBEI} for the special case $F (t)= e^{-\sum_{x\in
\Lambda}v_xt_x}$ when ${\rm Re}\, v_x \ge 0$.
Therefore, for this choice of $F$,
\begin{equation}
    \int_{[0,\infty)^{\carL }} F d\mu = \int_{[0,\infty)^{\carL }} F d\nu.
\end{equation}
This proves \eqref{e:DBEI} in the general case, since finite
measures are characterized by their Laplace transforms. The
hypothesis on the growth of $F$ assures its integrability.
\end{proof}

\subsection{Self-avoiding walk with loops}

Next, we derive a representation for a model of a self-avoiding walk
in a background of loops.  This requires the introduction of some
terminology and notation.

Given not necessarily distinct points $a,b \in \Lambda$,
  a \emph{self-avoiding walk} $\omega$
from $a$ to $b$ is a sequence
$x_0=a,x_1,x_2,\ldots, x_n=b$, for some $n \ge 1$, where $x_1,x_2,\ldots, x_{n-1}$ are distinct points in $\Lambda \setminus \{a,b \}$.  In other words,
for $a\not =b$,  $\omega$ is a non-intersecting path from $a$ to $b$ on the complete graph
on $\carL$ vertices and for $a=b$ it is non-intersecting except at $a=b$.  We again write $|\omega|$ for the length $n$ of $\omega$,
and sometimes regard $\omega$ as consisting of
directed edges rather than vertices.
Let $\Scal_{a,b}$
denote the set of all self-avoiding walks from $a$ to $b$.
For $X \subset \Lambda$,
we write $\Scal_{a,b}(X)$ for the subset of $\Scal_{a,b}$ consisting
of walks with $x_0=a$, $x_n=b$ and $x_1,x_2,\ldots,x_{n-1} \in X$.
A \emph{loop} $\gamma$ is an unrooted directed
cycle (consisting of distinct vertices)
in the complete graph, regarded sometimes as a cyclic list of vertices
and sometimes as directed edges.  We include the \emph{self-loop} which
joins a vertex to itself by a single edge,
as a possible loop (see Remark~\ref{rk:Wick} below).
We write $\Lcal$ for the set of all loops.
We write $\Gamma$ for a subgraph of $\Lambda$ consisting of
mutually-avoiding
loops, i.e., $\Gamma = \{ \gamma_1, \ldots, \gamma_m\}$
with each $\gamma_i \in \Lcal$ and $\gamma_i \cap \gamma_j =\varnothing$
(as sets of vertices) for $i \neq j$.
We write $\Gcal$ for the set of all such $\Gamma$
(including $\Gamma =\varnothing$),
and $\Gcal(X)$ for
the subset of $\Gcal$ which uses only vertices in $X \subset \Lambda$.
We write $|\gamma|$ for the length of $\gamma$, and
$|\Gamma|= \sum_{i=1}^m|\gamma_i|$ for the total length of loops
in~$\Gamma$.\looseness=-1

Given a $\Lambda \times \Lambda$ real matrix $C$,
$\omega \in \Wcal_{a,b}$ and $\Gamma \in \Gcal$, let
\begin{equation}
    C^\Gamma = \prod_{e\in\Gamma} C_{e},
    \quad\quad
    C^{\omega \cup \Gamma} = C^\omega C^\Gamma,
\end{equation}
where here we regard self-avoiding walks and loops as collections of directed edges
and write, e.g., $e = (\omega(i-1),\omega(i))$.
An empty product is equal to $1$.
We define the two-point function
\eqalign
\label{e:Gloop}
    \Gloop_{a,b}
    & =
    \sum_{\omega \in \Scal_{a,b}}
    \sum_{\Gamma \in \Gcal(\Lambda \setminus \omega)}
    C^{\omega \cup \Gamma}.
\enalign
The representation for $\Gloop_{a,b}$ is elementary and we derive it now.
\eject


\begin{theorem}
\label{thm:Gloop}
Let $C$ have positive Hermitian part.
Let $a,b \in \Lambda$ (not necessarily distinct) and let $X \subset \Lambda\setminus \{a,b\}$.  Then
\eqalign
\label{e:Sloop}
    \int_{\R^{2\carL}}
    d\mu_C  \bar{\phi}_a \phi_b
    \prod_{x \in X} (1+\phi_x\phib_x)
    &=
    \sum_{\omega \in \Scal_{a,b}(X )} C^\omega
    \int_{\R^{2\carL}}
    d\mu_C
    \prod_{x \in X \setminus \omega} (1+\phi_x\phib_x),
\\[-4pt]
\label{e:loops}
    \int_{\R^{2\carL}}
    d\mu_C
    \prod_{x \in X} (1+\phi_x\phib_x)
    &=
    \sum_{\Gamma \in \Gcal(X)} C^\Gamma,
\enalign
and, finally,
\begin{equation}
\label{e:GRloop}
    \Gloop_{a,b} =
    \int_{\R^{2\carL}}
    d\mu_C \phib_a \phi_b
    \prod_{x \in \Lambda: x \neq a,b} (1+\phi_x\phib_x).
\end{equation}
\end{theorem}

\begin{proof}
To prove \eqref{e:Sloop}, we write $F =\phi_b
\prod_{x \in X} (1+\phi_x\phib_x)$ and apply the integration
by parts formula \eqref{e:Gibp}, which replaces $\phib_aF$ by
$\sum_{v\in\Lambda}C_{a,v}\partial F/ \partial \phi_v$.  The first step in
the walk $\omega$ is $(a,v)$.
If the derivative acts on a factor in the product over $x$,
then it replaces that factor by $\phib_v$, and the procedure can
be iterated until the derivative acts on $\phi_b$, in which case
$\omega$ terminates.
The result is \eqref{e:Sloop}.

For \eqref{e:loops}, we expand the product to obtain
\begin{equation}
    \prod_{x \in X} (1+\phi_x\phib_x)
    =
    \sum_{Y \subset X} \prod_{y\in Y} \phi_y\phib_y.
\end{equation}
and hence
\begin{equation}
    \int_{\R^{2\carL}}
    d\mu_C
    \prod_{x \in X} (1+\phi_x\phib_x)
    =
    \sum_{Y \subset X}
    \int_{\R^{2\carL}}
    d\mu_C(u) \prod_{y\in Y} \phi_y\phib_y.
\end{equation}
We then evaluate the integral on the right-hand side using Lemma~\ref{lem:Wick},
and this gives \eqref{e:loops}.

The representation \eqref{e:GRloop} follows
from the combination of \eqref{e:Sloop}--\eqref{e:loops}.
\end{proof}

\begin{rk}
\label{rk:Wick}
Self-loops can be eliminated in the representation by replacing
the right-hand side of \eqref{e:GRloop} by
\begin{equation}
    \int_{\R^{2\carL}}
    d\mu_C \phib_a \phi_b \prod_{w \in \Lambda: x \neq a,b}
    (1+\wick{\phi_x\phib_x}),
\end{equation}
where
\begin{equation}
    \wick{\phi_x\phib_x} = \phi_x\phib_x - C_{x,x},
\end{equation}
using a modification of the above proof.
\end{rk}

\section{Self-avoiding walk models}
\label{sec:sawm}

\subsection{Self-avoiding walk}

We define the two-point function:
\eqalign
\label{e:Gsaw}
    \Gsaw_{a,b}
    & =
    \sum_{\omega \in \Scal_{a,b}}  C^\omega.
\enalign
When $a=b$, the walks are self-avoiding except for
the fact that the walk begins and ends at the same site.
In this case, there is, in particular,
a contribution due to the one-step walk that steps
from $a$ to $a$, which has weight $C_{a,a} \neq 0$.
The only new result in this paper is the integral representation for
$\Gsaw_{a,b}$.
The representation for the loop model \eqref{e:Gloop}
is easier than for \eqref{e:Gsaw}, as \eqref{e:Gloop} is in terms
of a bosonic (ordinary) Gaussian integral.
To eliminate the loops and obtain a representation
for the walk model \eqref{e:Gsaw}, we will need fermionic
(Grassmann) integrals involving anticommuting variables.
The necessary mathematical background for this is developed in
Section~\ref{sec:rwdf}, and the representation is stated and
derived in Section~\ref{sec:saw}.
This representation
is the point of departure for the analysis of the 4-dimensional self-avoiding
walk in \cite{BS10}, for a convenient particular choice of $C$.

\subsection{Weakly self-avoiding walk}

The two-point functions \eqref{e:Gloop} and \eqref{e:Gsaw}
are for strictly self-avoiding walks
and loops.  We also consider the
continuous-time weakly
self-avoiding walk, which is defined as follows.

Let $D$ have diagonal entries $d_x>0$,
$J$ have zero diagonal entries and $J_{x,y} \ge 0$,
and suppose that $D-J$ is diagonally
dominant.
Let $X$ and $\Ex_a$ be the continuous-time Markov process and corresponding
expectation, as in Theorem~\ref{thm:srw2}.  In particular, the process
dies at the random time $\zeta$ at which it makes a transition to the
cemetery state.
The local time at $x$ is given by
$L_{x} = \int_0^\infty \1_{X(s)=x} ds$ (note that the integral effectively
terminates at $\zeta < \infty$).
By definition,
\eqalign
    \sum_{x\in\Lambda} L_{x}^2
    &=
    \int_0^\infty ds_1  \int_0^\infty ds_2
    \sum_{x\in\Lambda} \1_{X(s_1)=x} \1_{X(s_2)=x}
    \\
   & =
    \int_0^\infty ds_1  \int_0^\infty ds_2
    \1_{X(s_1)= X(s_2)\neq \partial},
\enalign
so $\sum_{x\in\Lambda} L_{x}^2$ is a measure of the amount of self-intersection of
$X$ up to time $\zeta$.
The continuous-time weakly self-avoiding walk two-point function is
defined by
\eqalign
\label{e:Gwsaw}
    \Gwsaw_{a,b}
    & =
    \frac{1}{d_b\pi_{b,\partial}} \Ex_a
    \left(
    e^{-g\sum_{x\in \Lambda} L_{x}^{2}} e^{ - \lambda \zeta}
    \1_{X(\zeta^-) =b}
    \right),
\enalign
where $g>0$,
and $\lambda$ is a parameter (possibly negative) which is chosen
in such a way that the integral converges.
In \eqref{e:Gwsaw}, self-intersections
are suppressed by the factor
$\exp[-g\sum_{x\in \Lambda} L_{x}^{2}]$.
We will derive
a representation  for \eqref{e:Gwsaw} in Section~\ref{sec:wsaw}.

It follows from Proposition~\ref{prop:DBEI} that there is also the
alternate representation:
\begin{equation}
\label{e:GwsawDynkin}
    \Gwsaw_{a,b}
    =
    \int_0^\infty
    \overline\Ex_{a} \left(
    e^{-g\sum_{x\in \Lambda} L_{x,T}^{2}}
    e^{- \sum_x(\lambda +d_x - \overline d_x)L_{x,T}}
    \1_{X(T)=b} \right)
    dT.
\end{equation}
In the homogeneous case, in which $d_x-\overline d_x =a$ is independent
of $x$, the second exponential can be written as $e^{-\lambda'T}$
where $\lambda'=\lambda + a$.
This representation
is the starting point for the analysis of the weakly self-avoiding walk
on a 4-dimensional hierarchical lattice in \cite{BEI92,BI03c,BI03d},
on $\Z^4$ in \cite{BS10}, and for a model on $\Z^3$ in \cite{MS08}.

\section{Gaussian integrals with fermions}
\label{sec:rwdf}

In this section, we review some standard material about Gaussian
integrals which incorporate anticommuting Grassmann variables.  We
realize these Grassmann variables as differential forms.

\subsection{Differential forms}

We recall and extend the formalism introduced in Section~\ref{sec:wlm}.
Let $\Lambda = \{1,\ldots, \carL\}$ be a finite set of cardinality $\carL$.
Let $u_1,v_1,\ldots, u_\carL,v_\carL$ be standard coordinates on $\R^{2N}$,
so that $du_1 \wedge dv_1 \wedge \cdots \wedge du_\carL \wedge dv_\carL$ is the
standard volume form on $\R^{2\carL}$, where $\wedge$ denotes the usual
anticommuting wedge product (see \cite[Chapter~10]{Rudi76} for
an introduction).
We will drop the wedge from the notation
and write simply $du_idv_j$ in place of $du_i \wedge dv_j$.
The one-forms $du_i$, $dv_j$ generate the
Grassmann algebra of differential forms
on $\R^{2\carL}$.
A form which is a function of $u,v$ times a product of $p$ differentials
is said to have \emph{degree} $p$, for $p \ge 0$.

The integral of a differential form
over $\R^{2\carL}$ is defined to be zero unless the form has degree
$2\carL$.
A form $K$ of degree $2\carL$ can be written as
$K=f(u,v)du_1dv_1\cdots du_\carL dv_\carL$, and we define
\begin{equation}
    \int K = \int_{{\mathbb R}^{2\carL}} f(u,v)du_1dv_1\cdots du_\carL dv_\carL,
\end{equation}
where the right-hand side is the usual Lebesgue
integral of $f$ over $\R^{2\carL}$.

We again
complexify by setting $\phi_x = u_x + i v_x$, $\phib_x = u_x-iv_x$
and $d\phi_x = du_x+idv_x$, $d\phib_x = du_x-idv_x$, for $x \in \Lambda$.
Since the wedge product is
anticommutative, the following pairs all anticommute for every $x,y\in \Lambda$:
$d\phi_x$ and $d\phi_y$, $d\phib_x$ and $d\phi_y$,
$d\phib_x$ and $d\phib_y$.
Given an $\carL\times \carL$ matrix $A$, we write
$\phi A \phib = \sum_{x,y \in \Lambda} \phi_x A_{x,y} \phib_y$.
As in \eqref{e:ddphi},
\begin{equation}
\label{e:duv}
    d\phib_x d\phi_x = 2i du_x dv_x.
\end{equation}

The integral of a function $f(\phi,\phib)$ (a \emph{zero form})
with respect to $\prod_{x\in \Lambda}d\phib_xd\phi_x$ is thus given by
$(2i)^\carL$ times the integral of $f(u+iv,u-iv)$ over $\R^{2\carL}$.  Note that
the product over $x$ can be taken in any order, since each factor
$d\phib_xd\phi_x$
has even degree (namely degree two).
To simplify notation, it is convenient to introduce
\begin{equation}
\label{e:phipsi}
    \psi_x = \frac{1}{(2\pi i)^{1/2}} d \phi_x,
    \quad
    \psib_x = \frac{1}{(2\pi i)^{1/2}} d \phib_x,
\end{equation}
where we fix a choice of the square root and use this choice henceforth.
Then
\begin{equation}
    \psib_x \psi_x =\frac{1}{2\pi i} d\phib_x d\phi_x = \frac{1}{\pi}du_x dv_x.
\end{equation}

Given any matrix $A$, the \emph{action} is the even form defined by
\begin{equation}
\label{e:Sdef}
    S_A =  \phi A\bar{\phi}  + \psi A \bar{\psi} .
\end{equation}
In the special case $A_{u,v} = \delta_{u,x}\delta_{x,v}$, $S_A$ becomes
the form $\tau_x$ defined by
\begin{equation}
\label{e:taudef1}
    \tau_x = \phi_x \bar{\phi}_x + \psi_x \bar{\psi}_x.
\end{equation}

Let $K=(K_j)_{j\in J}$ be a collection of
forms.  When each $K_j$ is a sum of forms of even degree,
we say that $K$ is {\em even}.
Let $K_j^{(0)}$ denote the degree-zero part
of $K_j$.
Given a $C^\infty$ function $F : \R^{J} \to\mathbb C$
we define $F(K)$ by
its power series about the degree-zero part of $K$, i.e.,
\begin{equation}
\label{e:Fdef}
    F(K) = \sum_{\alpha} \frac{1}{\alpha !}
    F^{(\alpha)}(K^{(0)})
    (K - K^{(0)})^{\alpha}.
\end{equation}
Here $\alpha$ is a multi-index, with $\alpha ! = \prod_{j\in J}\alpha_j !$,
and
$(K - K^{(0)})^{\alpha}
=\prod_{j\in J} (K_{j} - K_{j}^{(0)})^{\alpha_{j}}$.
Note that the summation terminates as soon as
$\sum_{j\in J}\alpha_j=\carL$ since higher order forms vanish,
and that the order of the product on the right-hand side is irrelevant
when $K$ is even.
For example,
\begin{equation}
\label{e:eSA}
    e^{-S_A} = e^{-\phi A\phib}
    \sum_{n=0}^\carL \frac{(-1)^n}{ n!}
    (\psi A \psib)^n
    .
\end{equation}
Because the formal power series of a composition of two functions
is the same as the composition of the two formal power series,
we may regard $e^{-S_A}$ either as a function of the single
form $S_A$ or of the $\carL^2$
forms $\phi_x\phib_y + \frac{1}{2\pi i}d\phi_xd\phib_y$.
The same result is obtained for $e^{-S_A}$ in either
case.

\subsection{Gaussian integrals}

We refer to the integral $\int e^{-S_A} K$ as the mixed
bosonic-fermionic Gaussian expectation of $K$, or, more briefly, as
a mixed expectation. The following proposition shows that if $K$ is
a product of a zero form and factors of $\psi$ and $\psib$ then the
mixed expectation factorizes.  Moreover, if $K$ is a zero form then
the mixed expectation is just the usual Gaussian expectation of $K$,
and if $K$ is a product of factors of $\psi$ and $\psib$ then its
expectation is a determinant. It also shows that $\int e^{-S_A}$ is
self-normalizing in the sense that it is equal to $1$ without any
normalization required. The determinant in \eqref{e:Efac} appears
also e.g.\ in \cite[Lemma~B.7]{Salm99}, in a related purely
fermionic context and with a different proof.

\begin{prop}
\label{prop:sn}
Let $A$ have positive Hermitian part, with inverse $C=A^{-1}$.
Suppose that $f$ is a zero form.
Let $F=\prod_{r=1}^p\psib_{i_r}\prod_{s=1}^q \psi_{j_s}$.
If $p\neq q$ then $\int e^{-S_A} fF=0$.  When $p=q$, up to sign we
can take $F=\psib_{i_1}\psi_{j_1}\cdots\psib_{i_p}\psi_{j_p}$ and
in this case
\begin{equation}
\label{e:Efac}
    \int e^{-S_A} fF
    = \bigg( \int e^{-S_A}  f \bigg)
    \bigg( \int e^{-S_A} F \bigg)
    = I_f \det C_{i_1,\ldots,i_p;j_1,\ldots,j_p}
\end{equation}
where $I_f=\int f d\mu_C(\phi,\phib)$,
and where $C_{i_1,\ldots,i_p;j_1,\ldots,j_p}$
is the $p\times p$ matrix whose $r,s$ element
is $C_{i_r,j_s}$ when $p\neq 0$,
and the determinant is replaced by $1$ when $p=0$.
In particular,
\begin{equation}
\label{e:sn}
    \int e^{-S_A}=1.
\end{equation}
\end{prop}

\begin{proof}
We first note that if $p \neq q$ then
no form of degree $2\carL$ can be obtained by expanding $e^{- \psi
A \bar\psi}F$ and the integral vanishes.  Thus we assume $p=q$.

Let $i = i_1,\ldots, i_p$, $j = j_1,\ldots, j_p$, and
\begin{equation}
    B_{i,j}
    =
    \int e^{-S_A}f \psib_{i_1}\psi_{j_1}\cdots\psib_{i_p}\psi_{j_p}.
\end{equation}
For $k \in \Lambda$, let
\begin{equation}
\label{e:psitildef}
    \tilde{\psi}_{k}
    =
    \sum_{l\in \Lambda } A_{k,l}\bar{\psi}_{l}.
\end{equation}
The tensor product $A^{\otimes p}$ is
a linear operator on $V^{\otimes p}$ defined by the matrix
elements
\begin{equation}
    (A^{\otimes p})_{i,j}
    =
    A_{i_{1},j_{1}}
    A_{i_{2},j_{2}}\dotsb
    A_{i_{p},j_{p}}.
\end{equation}
By definition, \eqref{e:eSA}, and the anticommutation relation
$\psi_{k_l} \tilde\psi_{k_l} = - \tilde\psi_{k_l} \psi_{k_l}$,
\eqalign
    &(A^{\otimes p}B)_{i,j}
     =
    \int e^{-S_A}f
    \tilde\psi_{i_1}\psi_{j_1}\cdots\tilde\psi_{i_p}\psi_{j_p}
    \nnb
    & =
    \frac{1}{(\carL -p)!} \sum_{k_1,\ldots,k_{\carL-p}}
    \int e^{-\phi A\phib}f
    \tilde\psi_{k_1}\psi_{k_1}\cdots\tilde\psi_{k_{\carL-p}}\psi_{k_{\carL-p}}
    \tilde\psi_{i_1}\psi_{j_1}\cdots\tilde\psi_{i_p}\psi_{j_p}.
\label{e:AB1}
\enalign
By antisymmetry, for a nonzero contribution,
$k_1,\ldots,k_{\carL-p},i_1,\ldots,i_p$ must be
a permutation of $\Lambda$, as must be
$k_1,\ldots,k_{\carL-p},j_1,\ldots,j_p$.  In particular,
$j_1,\ldots,j_p$ must be a permutation of $i_1,\ldots , i_p$;
let $\epsilon_{i,j}$ be the sign of this permutation
(and equal zero if it is not a permutation).  Then we
can rearrange the above to obtain
\eqalign
\label{e:AB2}
    (A^{\otimes p}B)_{i,j}
    & =
    \epsilon_{i,j}\int e^{-\phi A\phib}f
    \tilde\psi_{1}\psi_{1}\cdots\tilde\psi_{\carL}\psi_{\carL}.
\enalign
We insert \eqref{e:psitildef} on the right-hand side and again use
antisymmetry and then Lemma~\ref{lem:ZC} to obtain
\eqalign
\label{e:AB3}
    (A^{\otimes p}B)_{i,j}
    & =
    \epsilon_{i,j}\det A
    \int e^{-\phi A\phib}f
    \psib_{1}\psi_{1}\cdots \psib_{\carL}\psi_{\carL}
    =
    I_f \epsilon_{i,j}  .
\enalign
When $p=0$ the above calculations give $B=I_f$, as required.

For $p \neq 0$, we use the fact that $C^{\otimes p}$ is the inverse
of $A^{\otimes p}$ to obtain
\begin{equation}
    B_{k,j} = \sum_{l}C^{\otimes p}_{k,l} (A^{\otimes p} B)_{l,j}
    = I_f \sum_{l} \epsilon_{l,j} C^{\otimes p}_{k,l}.
\end{equation}
The sum on the right-hand side is the determinant
$\det C_{k_1,\ldots,k_p;j_1,\ldots,j_p}$, as required.
\end{proof}

%

In the Gaussian integral in the above proposition, the fermionic
part $d\phi A d\phib$ of the action gives rise to a factor $\det A$
while the bosonic part $\phi A \phib$ gives rise to the reciprocal
of this determinant, providing the cancellation that produces the
self-normalization property \eqref{e:sn}.

We will use the following corollary in Section~\ref{sec:expres}.

\begin{cor}
\label{cor:Gram1ii}
Let $x_1,\ldots, x_k$ be distinct elements of $\Lambda$.
Then
\begin{equation}
\label{e:Gramloop}
     \int e^{-S_A} \psi_{x_1}\psib_{x_1}\cdots \psi_{x_k}\psib_{x_k}
     =
     \sum_{\sigma \in S_k} (-1)^{N(\sigma)} \prod_{l=1}^k C_{x_l,\sigma(x_l)},
\end{equation}
where $N(\sigma)$ is the number of cycles in the permutation $\sigma$.
\end{cor}

\proof
It follows from \eqref{e:Efac} and anticommutativity that
\begin{equation}
     \int e^{-S_A} \psi_{x_1}\psib_{x_1}\cdots \psi_{x_k}\psib_{x_k}
     =
     (-1)^k
     \sum_{\sigma \in S_k} \epsilon_\sigma \prod_{l=1}^k C_{x_l,\sigma(x_l)},
\end{equation}
where $\epsilon_\sigma$ is the sign of the permutation $\sigma$.
Then \eqref{e:Gramloop} follows from the identity
\begin{equation}
    \epsilon_\sigma = (-1)^k (-1)^{N(\sigma)},
\end{equation}
which itself follows from the fact that for a permutation $\sigma \in S_k$
consisting of cycles $c$ of length $|c|$,
\begin{equation}
    \epsilon_\sigma = \prod_{c \in \sigma} \epsilon_{c}
    =
    \prod_{c\in \sigma} (-1)^{|c|+1} = (-1)^k (-1)^{N(\sigma)}.
\end{equation}
\qed

\begin{rk}
The omission of the operation $A^{\otimes p}$ in \eqref{e:AB1}--\eqref{e:AB3}
leads to the alternate formula
\begin{equation}
\label{e:Bijalt}
    B_{i,j}
    =
    \int e^{-S_A}f \psib_{i_1}\psi_{j_1}\cdots\psib_{i_p}\psi_{j_p}
    =
    I_f
    \frac{1}{\det A}  \det \hat A_{i_1,\ldots,i_p;j_1,\ldots,j_p}
    \epsilon_{\sigma_i}\epsilon_{\sigma_j},
\end{equation}
where $\sigma_i \in S_\carL$ is the permutation that moves $i_1,\ldots
, i_p$ to $1,\ldots,p$ and preserves the order of the other indices
and $\epsilon_{\sigma_i}$ is its sign
(and similarly for $\sigma_j$), and
where $\hat A_{i_1,\ldots,i_p;j_1,\ldots,j_p}$ is the
$(\carL-p)\times (\carL-p)$ matrix obtained
from
$A$ by deleting rows $j_1,\ldots, j_p$ and columns $i_1,\ldots, i_p$.
The identity \eqref{e:Bijalt} is essentially \cite[Lemma~4]{BM91}.
This proves the
fact from linear algebra that
\begin{equation}
\label{e:Cramer}
    \det C_{i_1,\ldots,i_p;j_1,\ldots,j_p}
    =
    \frac{1}{\det A}  \det \hat A_{i_1,\ldots,i_p;j_1,\ldots,j_p}
    \epsilon_{\sigma_i}\epsilon_{\sigma_j}.
\end{equation}
The case $p=1$ of \eqref{e:Cramer} states that
\begin{equation}
\label{e:Cramer1}
    C_{i_1;j_1} = A^{-1}_{i_1;j_1}
    = \frac{1}{\det A}  \det \hat A_{i_1;j_1}(-1)^{i_1+j_1},
\end{equation}
which is Cramer's rule.  Thus \eqref{e:Cramer} is a generalization
of Cramer's rule.
\end{rk}

\subsection{Integrals of functions of $\tau$}

The identity \eqref{e:F0} below provides an extension of
\eqref{e:sn}, and will be used in Section~\ref{sec:saw}.
The identity \eqref{e:FDynkin} is sometimes called the $\tau$-isomorphism;
it will lead to a representation for the weakly self-avoiding
walk two-point function \eqref{e:Gwsaw}.
Our method of proof follows the method of \cite{BEI92,Imbr03}.
Alternate approaches to \eqref{e:F0} are given in Sections~\ref{sec:expres}
and \ref{sec:susy}.

Recall the definitions of $\tau_x$ in \eqref{e:taudef1}
and $L_x$ above Theorem~\ref{thm:srw2}.
We write $\tau$ for the entire collection $(\tau_x)_{x\in\Lambda}$,
and similarly for $L$.


\begin{prop}
\label{prop:Ftau0}
Suppose that $A$ has positive Hermitian part.
Let $F$ be a $C^\infty$ function on $[0,\infty)^\carL$
($C^\infty$ also on the boundary), and assume that for each
$\epsilon >0$  and  multi-index $\alpha$
there is a constant $C=C_{\epsilon,\alpha}$ such that
$F$ and its derivatives obey
$|F^{(\alpha)} (t)| \le C
\exp (\epsilon \sum_{x\in\Lambda}t_{x})$ for all
$t \in [0,\infty)^\carL$.
Then
\eqalign
\label{e:F0}
    \int e^{-S_A} F(\tau) &= F(0).
\enalign
Suppose further that $A=D-J$ is diagonally dominant and real.
Then
\eqalign
\label{e:FDynkin}
    \int e^{-S_A} F(\tau) \bar\phi_a \phi_b
    &=
    \frac{1}{d_b \pi_{b,\partial}}
    \Ex_a \left( F(L) \1_{X(\zeta^-) =b}\right).
\enalign
\end{prop}

\proof It is straightforward to adapt the result of \cite{Seel64} to
extend $F$ to a $C^\infty$ function on $\R^\carL$, which we also
call $F$. By multiplying $F$ by a suitable $C^\infty$ function, we
can further assume that $F$ is equal to zero on the complement of
$[-1,\infty)^\carL$. Fix $\epsilon >0$ such that $A-\epsilon I$ has
positive Hermitian part, and let $H (t) = F (t)\exp (-\epsilon
\sum_{x}t_{x})$. Then $H$ is a Schwartz class function.  Its Fourier
transform is defined by
\begin{equation}
    \widehat H (v) = \int_{\R^\carL} H(t) e^{i v \cdot t}
    dt_1 \ldots dt_\carL,
\end{equation}
where $v \cdot t =\sum_{x\in \Lambda}v_xt_x$.
The function $H$ can be recovered via the inverse Fourier transform
as
\begin{equation}
    H (t)
    =
    (2\pi )^{-M} \int_{\R^\carL} \widehat{H} (v)e^{-iv \cdot t}
    \,dv_1 \ldots dv_\carL.
\end{equation}
Since $H$ is of Schwartz class, the above integral is absolutely
convergent.  Also,
\begin{equation}
\label{e:invFT}
    F (t)
    =
    (2\pi )^{-M} \int_{\R^\carL}
    \widehat{H} (v)e^{\sum_{x} (-iv_{x}+\epsilon )t_{x}}
    \,dv_1 \ldots dv_\carL .
\end{equation}

We may replace $t$ by $\tau$ in \eqref{e:invFT}; this amounts to a
statement about differentiating under the integral since functions
of $\tau$ are defined by their power series as in \eqref{e:Fdef}.
Let $V$ be the real diagonal matrix with $V_{x,x}=v_x$.  Since
$A-\epsilon I +iV$ has positive Hermitian part, \eqref{e:sn} gives
\begin{equation}
    \int e^{-S_A} e^{\sum_{x} (-iv_{x}+\epsilon )\tau_{x}}
    = \int e^{-S_{A-\epsilon I +iV}} = 1.
\end{equation}
Assuming that it is possible to interchange the integrals, we obtain
\begin{equation}
    \int e^{-S_A} F(\tau)
    =
    (2\pi )^{-M} \int_{\R^\carL}
    \widehat{H} (v)
    dv_1 \ldots dv_\carL
    = H(0)=F(0),
\end{equation}
which is \eqref{e:F0}.

To complete the proof of \eqref{e:F0},
it remains only to justify the interchange
of integrals; this can be done as follows.
By definition, the iterated integral
\begin{equation}
\label{e:fub}
    \int e^{-S_A}
    \int_{\R^\carL}
    dv_1 \ldots dv_\carL\,
    \widehat{H} (v)e^{\sum_{x} (-iv_{x}+\epsilon )\tau_{x}}
\end{equation}
is equal to
\eqalign
    \sum_{n,N} \frac{(-1)^N}{n!N!} & \int e^{-\phi A \phib}
    (\psi A \psib)^N
    \left(\sum_x (-iv_{x}+\epsilon )\psi_{x}\psib_x\right)^n \nnb
    &\times  \int_{\R^\carL}
    dv_1 \ldots dv_\carL\,
    \widehat{H} (v)e^{\sum_{x} (-iv_{x}+\epsilon )\phi_x\phib_{x}}.
\enalign
According to our definition of integration, the outer integral is
evaluated as a usual Lebesgue integral by keeping the (finitely
many) terms that produce the standard volume form on $\R^{2\carL}$.
Since $\widehat{H}$ is Schwartz class and $A-\epsilon I$ has
positive Hermitian part, the resulting iterated Lebesgue integral is
absolutely convergent and its order can be interchanged by Fubini's
theorem.  Once the integrals have been interchanged, the sums over
$n$ and $N$ can be resummed to see that \eqref{e:fub} has the same
value when its two integrals are interchanged, and the proof of
\eqref{e:F0} is complete.

To prove \eqref{e:FDynkin}, we fix $\epsilon >0$ such that
$A-\epsilon I$ is diagonally dominant.
Then
\eqalign
    \int e^{-S_A} e^{\sum_{x} (-iv_{x}+\epsilon )\tau_{x}}
    \bar\phi_a \phi_b
    & =
    \int e^{-S_{A-\epsilon I + iV}} \bar\phi_a \phi_b
    =
    \Gsrw_{a,b}(-\epsilon + iv)
    \nnb
    & =
    \frac{1}{d_b\pi_{b,\partial}}
    \Ex_{a}\left( e^{(\epsilon -iv)\cdot L} \1_{X(\zeta^-)=b} \right)
    ,
\enalign
where we have used \eqref{e:Efac} and Theorem~\ref{thm:srw}
in the second equality, and
Theorem~\ref{thm:srw2} in the third.
With further application of Fubini's theorem, we obtain
\eqalign
    \int e^{-S_A} \bar\phi_a \phi_b F(\tau)
    & =
    \frac{1}{d_b\pi_{b,\partial}}\Ex_{a} \left(
    e^{\epsilon \cdot L}
    (2\pi )^{-M} \int_{\R^\carL}
    \widehat{H} (v)
    e^{-iv \cdot L}dv \; \1_{X(\zeta^-)=b}
    \right)
    \nnb
    & =
    \frac{1}{d_b\pi_{b,\partial}}
    \Ex_{a} \left( F(L)  \1_{X(\zeta^-)=b} \right) ,
\enalign
which is \eqref{e:FDynkin}.
\qed

\section{Self-avoiding walk representations}
\label{sec:fermion}

\subsection{Weakly self-avoiding walk}
\label{sec:wsaw}

\subsubsection{The representation}

\begin{theorem}
\label{thm:wsaw}
The weakly self-avoiding walk two-point function
$\Gwsaw_{a,b}$ has the
representation
\begin{equation}
\label{e:Gsawrep}
    \Gwsaw_{a,b}
    =
    \int e^{-S_A} \bar\phi_a \phi_b
    e^{-g\sum_{x\in\Lambda}\tau_x^2-\lambda \sum_{x\in\Lambda}\tau_x}.
\end{equation}
\end{theorem}

\begin{proof}
This is immediate when we
take
$F(\tau)=e^{-g\sum_{x\in\Lambda}\tau_x^2-\lambda \sum_{x\in\Lambda}\tau_x}$
in \eqref{e:FDynkin}, and compare with \eqref{e:Gwsaw}.
\end{proof}

\subsubsection{The $N \to 0$ limit}
\label{sec:N0lim}

If we omit the fermions from the right-hand side of
\eqref{e:Gsawrep} and normalize the integral then we obtain instead
the two-point function of the $|\phi|^4$ field theory, namely
\begin{equation}
    \langle \phib_a \phi_b \rangle
    =
    \frac{
    \int d\mu_C \phib_a \phi_b
    e^{-g\sum_{x\in\Lambda}|\phi_x|^4-\lambda \sum_{x\in\Lambda}|\phi_x|^2}
    }
    {
    \int d\mu_C
    e^{-g\sum_{x\in\Lambda}|\phi_x|^4-\lambda \sum_{x\in\Lambda}|\phi_x|^2}
    }.
\end{equation}
This is known to have a representation as the two-point function of
a system of a weakly self-avoiding walk and weakly self-avoiding
loops, all weakly mutually-avoiding \cite{BFS82,Syma69}, as we now
briefly sketch.

Let $n_x(\omega)$ denote the number of visits to $x$ by a
walk $\omega$.
Let
\begin{equation}
\label{e:dnudef}
    d\nu_n(s) =
    \begin{cases}
    \delta(s) ds & n=0 \\[3pt]
    \frac{s^{n-1}}{(n-1)!} \1_{s \ge 0} ds & n \ge 1
    \end{cases},
    \quad\quad\quad
    d\nu_\omega(t)
    =
    \prod_{x \in \Lambda} d\nu_{n_x(\omega)}(t_x).
\end{equation}
It follows from \cite[Theorem~2.1]{BFS82}
(see also \cite[p.137]{BFS83II} and \cite[p.197]{FFS92})
that for a real $N$-component field $\phi$, for any component $i$ we have
\eqalign
   \langle \phi_a^{(i)}\phi_b^{(i)}
    \rangle
    &=
    \frac{1}{Z}
    \sum_{n=0}^\infty \frac{1}{n!} \left( \frac N2\right)^n
    \sum_{\omega \in \Wcal_{a,b}}
    \sum_{x_1,\ldots,x_n\in \Lambda}
    \sum_{\omega_1\in \Wcal_{x_1,x_1}}\cdots \sum_{\omega_n \in \Wcal_{x_n,x_n}}
    \nnb& \times \frac{J^{\omega \cup \omega_1 \cup \cdots \cup \omega_n}}
    {\|\omega_1\| \cdots \|\omega_n\|}
    \int d\nu_{\omega \cup \omega_1 \cup \cdots \cup \omega_n}(t)
    e^{-4g \sum_{x\in \Lambda} t_x^2 - 2\lambda \sum_{x\in \Lambda}t_x},
\label{e:phi4}
\enalign
where $\|\omega\|=|\omega|+1$ denotes the number of vertices in $\omega$,
\eqalign
    Z
    &
    = \sum_{n=0}^\infty \frac{1}{n!} \left( \frac N2\right)^n
    \sum_{x_1,\ldots,x_n\in \Lambda}
    \sum_{\omega_1\in \Wcal_{x_1,x_1}}\cdots \sum_{\omega_n \in \Wcal_{x_n,x_n}}
   \nnb &  \times
    \frac{J^{ \omega_1 \cup \cdots \cup \omega_n}}
    {\|\omega_1\| \cdots \|\omega_n\|}
    \int d\nu_{\omega_1 \cup \cdots \cup \omega_n}(t)
    e^{-4g \sum_{x\in \Lambda} t_x^2 - 2\lambda \sum_{x\in \Lambda}t_x}
\enalign is a normalization constant, and
\begin{equation}
    d\nu_{\omega \cup \omega_1 \cup \cdots \cup \omega_n}(t)
    =
    \prod_{x \in \Lambda}
    d\nu_{n_x(\omega) + n_x(\omega_1) + \cdots + n_x(\omega_n)+N/2}(t_x).
\end{equation}

Note the factor $N/2$ associated to each loop.
If we simply set $N=0$ in these formulas, then only the $n=0$ term survives,
and we obtain the formal limit (formal, because the left-hand side is
defined only for $N=1,2,3,\ldots$)
\begin{equation}
\label{e:N0limit}
    \lim_{N \to 0}
    \langle \phi_a^{(1)}\phi_b^{(1)}
    \rangle
    =
    \sum_{\omega \in \Wcal_{a,b}}
    J^{\omega}
    \int
    d\nu_{\omega}(t)
    e^{-4g \sum_{x\in \Lambda} t_x^2 - 2\lambda \sum_{x\in \Lambda}t_x}.
\end{equation}
As we argue next, the right-hand side of \eqref{e:N0limit}
is equal to the weakly self-avoiding walk two-point function
$\Gwsaw_{a,b}$ (with modified parameters $g,\lambda$).
This recovers de Gennes' idea, in the context of the weakly self-avoiding
walk \cite{ACF83}.

We now
show that the right-hand side of \eqref{e:N0limit} is equal
to the right-hand side in the representation \eqref{e:Gwsaw}
of $\Gwsaw_{a,b}$, with constant $d_x\equiv d$.  As in
the proof of Theorem~\ref{thm:srw2}, we condition on the events
$\{\eta =n\}$ and also on $Y=(Y_0,Y_1,\ldots, Y_n)\in \Wcal_{a,b}$.
Given both of these,
the random variable $L_{x}$ has a $\Gamma(n_x(Y),d)$ distribution,
since it is the sum of independent ${\rm Exp}(d)$ random variables.
Thus we obtain
\eqalign
    \Gwsaw_{a,b}
    & =
    \frac{1}{d\pi_{b,\partial}} \Ex_a
    \left(
    e^{-g \sum_x L_x^2 - \lambda \sum_{x}L_x}
    \1_{X(\zeta^-) =b}
    \right)
    \nnb
    & =
    \frac{1}{d} \sum_{n=0}^\infty
    \Ex_a\left[\Ex_a
    \left(
    e^{-g \sum_x L_{x}^2 - \lambda \sum_{x}L_{x}}
    | Y_0,\ldots, Y_n
    \right) \right].
\enalign
Since
\eqalign
    \Ex_a
    \left(
    e^{-g \sum_x L_{x}^2 - \lambda \sum_{x}L_{x}}
    | Y_0,\ldots, Y_n
    \right)
    & =
    \int d\Gamma_Y (t) e^{-g\sum_x t_x^2 - \lambda\sum_x t_x}
\enalign
with
\begin{equation}
    d\Gamma_Y (t)
    =
    \prod_{x \in \Lambda} d\nu_{n_x(Y)}(t_x)
    d^{n_x(Y)} e^{-dt_x}
    =
    d\nu_Y(t) d^{n+1} e^{-d\sum_x t_x},
\end{equation}
this gives
\eqalign
    \Gwsaw_{a,b}
    & = \frac{1}{d} \sum_{n=0}^\infty
    \sum_{\omega \in \Wcal_{a,b} : |\omega|=n}
    \left( \frac{J}{d}\right)^\omega
    \int d\nu_\omega(t) d^{n+1} e^{-d\sum_x t_x}
    e^{-g\sum_x t_x^2 - \lambda\sum_x t_x}
    \nnb
    & =
    \sum_{\omega \in \Wcal_{a,b}}
    J^\omega
    \int d\nu_\omega(t)
    e^{-g\sum_x t_x^2 - (\lambda+d)\sum_x t_x},
\enalign
which is the right-hand side of \eqref{e:N0limit} with a modified
choice of constants in the exponent.

%

Theorem~\ref{thm:wsaw} provides an alternative to
the above formal $N \to 0$ limit.
The inclusion of fermions in Theorem~\ref{thm:wsaw}
has eliminated all the loops, leaving
only the weakly self-avoiding walk.
In Section~\ref{sec:expres}, we will\vadjust{\eject}
make explicit the mechanism by which this occurs in the
strictly self-avoiding walk
representation:
fermionic loops cancel the bosonic ones.

\subsection{Strictly self-avoiding walk}
\label{sec:saw}

Here we obtain the representation for \eqref{e:Gsaw}.
We give two proofs based on two different ideas.

\subsubsection{Proof by expansion and resummation}
\label{sec:expres}

\begin{theorem}
\label{thm:saw}
Let $A$ have positive Hermitian part,
and let $C=A^{-1}$ denote its inverse.
For all $a,b \in \Lambda$,
\begin{equation}
\label{e:Gsawid}
    \Gsaw_{a,b}
    =
    \int e^{-S_A} \bar{\phi}_a \phi_b
    \prod_{x \in \Lambda\setminus\{a,b\}}(1+\tau_x).
\end{equation}
\end{theorem}

\proof
We write $X = \Lambda \setminus \{a,b\}$.
By expanding the product of $1+\tau_x = (1+\phi_x\phib_x)+\psi_x\psib_x$,
we obtain
\begin{equation}
\label{e:tauprod}
    \prod_{x \in X} (1+\tau_x)
    =
    \sum_{Y \subset X} \left( \prod_{y \in Y} \psi_y\psib_y \right)
    \left( \prod_{z \in X \setminus Y} (1+\phi_z\phib_z)\right).
\end{equation}
Thus, by Proposition~\ref{prop:sn},
\eqalign
    &\int e^{-S_A} \bar{\phi}_a \phi_b
    \prod_{x \in X}(1+\tau_x) \nnb
    &\quad\quad=
    \sum_{Y \subset X} \left( \int e^{-S_A}\prod_{y \in Y} \psi_y\psib_y \right)
    \left(\int e^{-S_A} \bar{\phi}_a \phi_b
    \prod_{z \in X \setminus Y} (1+\phi_z\phib_z)\right).
\enalign
By \eqref{e:Sloop},
\begin{equation}
    \int e^{-S_A} \bar{\phi}_a \phi_b
    \prod_{z \in X \setminus Y} (1+\phi_z\phib_z)
    =
    \sum_{\omega \in \Scal_{a,b}(X \setminus Y)} C^\omega
    \int e^{-S_A}
    \prod_{z \in X \setminus (Y\cup\omega)} (1+\phi_z\phib_z),
\end{equation}
where we have also used \eqref{e:Efac} twice to equate bosonic and mixed
bosonic-fermionic integrals.
Another application of Proposition~\ref{prop:sn} then gives
\eqalign
   & \int e^{-S_A} \bar{\phi}_a \phi_b
    \prod_{x \in X}(1+\tau_x) \nnb
    &\quad\quad=
    \sum_{Y \subset X} \sum_{\omega \in \Scal_{a,b}(X \setminus Y)} C^\omega
    \int e^{-S_A}
    \prod_{y \in Y} \psi_y\psib_y
    \prod_{z \in X \setminus (Y\cup\omega)} (1+\phi_z\phib_z).
\enalign
We now interchange the sums over $Y$ and $\omega$, and then resum to obtain
\eqalign
   & \int e^{-S_A} \bar{\phi}_a \phi_b
    \prod_{x \in X}(1+\tau_x) \nnb[-2pt]
    &\quad\quad =
    \sum_{\omega \in \Scal_{a,b}} C^\omega
    \sum_{Y \subset X\setminus \omega}
    \int e^{-S_A}
    \prod_{y \in Y} \psi_y\psib_y
    \prod_{z \in (X \setminus\omega)\setminus  Y} (1+\phi_z\phib_z)
    \nnb
    & \quad\quad=
    \sum_{\omega \in \Scal_{a,b}} C^\omega
    \int e^{-S_A}
    \prod_{x \in X \setminus\omega} (1+\tau_x).
\enalign
By \eqref{e:F0}, the integral in the last line is $1$, and we obtain
\eqref{e:Gsawid}.
\qed


The above proof ultimately relies on the identity
\begin{equation}
\label{e:loopcan}
    \int e^{-S_A} \prod_{x \in X} (1+\tau_x)
    =1,
\end{equation}
for a subset $X \subset \Lambda$.  This identity follows immediately
from \eqref{e:F0}.  We now give an alternate, more direct proof of
\eqref{e:loopcan}, which demonstrates that \eqref{e:loopcan} results
from the explicit cancellation of bosonic loops carrying a factor
$+1$ with fermionic loops carrying a factor $(-1)$.  The net effect
of a loop is $(+1)+(-1)=0$, which provides a realization of the
self-avoiding walk as corresponding to an $N=0$ model, without the
need of a mysterious $N \to 0$ limit.

\medskip \noindent
\emph{Alternate proof of \eqref{e:loopcan}.}
We expand the last product in \eqref{e:tauprod} and apply Proposition~\ref{prop:sn}
to obtain
\begin{equation}
\label{e:loopcan1}
    \int e^{-S_A} \prod_{x \in X}(1+\tau_x)
    =
    \sum_{{\rm disjoint} \, X_1,X_2 \subset X}
    \int e^{-S_A}
    \prod_{u \in X_1} \phi_u\phib_u
    \int e^{-S_A}
    \prod_{v \in X_2} \psi_v\psib_v.
\end{equation}
The term $X_1=X_2=\varnothing$ is special, and contributes $1$ to the
above right-hand side.
We write $S(X_i)$ for the set of permutations of $X_i$, $c_i$ for
a cycle of $\sigma_i \in S(X_i)$, and $W_{c_i} = \prod_{e \in c_i}C_e$
for the weight of the loop
corresponding to the cycle $c_i$.
With this notation, we can evaluate the integrals using
Lemma~\ref{lem:Wick} and \eqref{e:Gramloop} to find that the contribution to
the right-hand side of \eqref{e:loopcan1} due to all terms other than
$X_1=X_2=\varnothing$ is equal to
\begin{equation}
    \sum_{Y \subset X: Y \neq \varnothing}
    \sumtwo{{\rm disjoint} \, X_1,X_2 }{X_1\cup X_2 = Y}
    \sumtwo{\sigma_1 \in S(X_1)}{\sigma_2 \in S(X_2)}
    \prod_{c_1\in \sigma_1} W_{c_1}  \prod_{c_2\in \sigma_2} (-W_{c_2})
    .
\end{equation}
We claim that this equals
\begin{equation}
     \sum_{Y \subset X: Y \neq \varnothing}
    \sum_{\sigma \in S(Y)}
    \prod_{c\in \sigma} \left(
    W_{c}  + (- W_{c})    \right) = 0.
\end{equation}
This is a consequence of the fact that, for fixed $Y$,
\begin{equation}
    \sum_{\sigma \in S(Y)}
    \prod_{c\in \sigma} \left(
    P_{c}  + Q_{c}    \right)
    =\!\!
    \sumtwo{{\rm disjoint} \, X_1,X_2}{X_1 \cup X_2=Y}\!\!
    \sumtwo{\sigma_1 \in S(X_1)}{\sigma_2 \in S(X_2)}
    \prod_{c_1\in \sigma_1} P_{c_1}  \prod_{c_2\in \sigma_2} Q_{c_2}
    ,
\end{equation}
which follows by expanding the product on the left-hand
side.
\qed
\eject

\subsubsection{Proof by integration by parts}
\label{sec:ibp}

The integration by parts formula \eqref{e:Gibp} extends easily
to the mixed bosonic-fermionic case, to give
\begin{equation}
\label{e:ibp}
    \int e^{-S_A} \bar{\phi}_x F
    = \sum_{v\in\Lambda} C_{x,v} \int e^{-S_A}
    \frac{\partial F}{\partial \phi_v},
\end{equation}
where $A$ has positive Hermitian part, $C=A^{-1}$, and
$F$ is any $C^\infty$ form such that both sides are integrable.
To see this, we first note that by linearity it suffices to consider
the case $F=fK$ where $f$ is a zero form and $K$ is a product of
factors of $\psi$ and $\psib$.  By Proposition~\ref{prop:sn} and
\eqref{e:Gibp},
\eqalign
    \int e^{-S_A} \bar{\phi}_x fK
    &=
    \int e^{-S_A} \bar{\phi}_x f \int e^{-S_A} K
    \nnb
    &=
    \sum_{v\in\Lambda} C_{x,v} \int e^{-S_A}
    \frac{\partial f}{\partial \phi_v} \int e^{-S_A} K
    \nnb
    &=
    \sum_{v\in\Lambda} C_{x,v} \int e^{-S_A}
    \frac{\partial fK}{\partial \phi_v},
\enalign
and this proves \eqref{e:ibp}.

The special case $F=\phi_y$ in \eqref{e:ibp} gives
$\int e^{-S_A} \bar{\phi}_a \phi_b = C_{a,b}$.  More interestingly,
the choice $F= \phi_b(1+\tau_x)$ gives
$\int e^{-S_A} \bar{\phi}_a \phi_b  (1+\tau_x)= C_{a,b} +  C_{a,x}C_{x,b}$.
In the Gaussian integral, the effect of $\bar\phi_a$ is to start a walk
step at $a$, whereas $\phi_b$ has the effect of terminating a walk step
at $b$.  Each step receives the appropriate matrix element of the
covariance $C$ as its weight.  This leads to the following alternate
proof of Theorem~\ref{thm:saw}.


\smallskip \noindent
\emph{Second proof of Theorem~\ref{thm:saw}.}
The right-hand side of \eqref{e:Gsawid} is equal to
\begin{equation}
\label{e:zFint}
    \int e^{-S_A} \bar{\phi}_a F
\end{equation}
with
\begin{equation}
    F = \phi_b \prod_{x \neq a,b}(1+\tau_x),
\end{equation}
and hence
\begin{equation}
\label{e:Fder}
    \frac{\partial F}{\partial \phi_v}
    =
    \delta_{b,v} \prod_{x \neq a,b}(1+\tau_x)
    +
    \1_{v \ne a,b} \phi_b \bar{\phi}_v \prod_{x \neq a,b,v}(1+\tau_x).
\end{equation}
Substitution of \eqref{e:Fder} into \eqref{e:ibp}, using \eqref{e:F0},
gives
\begin{equation}
    \int e^{-S_A} \bar{\phi}_a F
    =
    C_{a,b} + \sum_{v\ne a,b} C_{a,v} \int e^{-S_A}
    \bar{\phi}_v \phi_b \prod_{x \neq a,b,v}(1+\tau_x).
\end{equation}
After iteration, the right-hand side gives $\Gsaw_{a,b}$.
\qed

\subsection{Comparison of two self-avoiding walk
representations}

The representations \eqref{e:Gsawrep} and \eqref{e:Gsawid} state that
\eqalign
\label{e:Gsawrep2}
    \Gwsaw_{a,b}
    &=
    \int e^{-S_A} \bar\phi_a \phi_b
    e^{-g\sum_{x\in\Lambda}\tau_x^2-\lambda \sum_{x\in\Lambda}\tau_x}.
\\
\label{e:Gsawid2}
    \Gsaw_{a,b}
    &=
    \int e^{-S_A} \bar{\phi}_a \phi_b
    \prod_{x \in \Lambda\setminus\{a,b\}}(1+\tau_x).
\enalign
These are heuristically related as follows.
We insert the missing factors for $x=a,b$ in the product
in \eqref{e:Gsawid2}, and make the (uncontrolled) approximation
\begin{equation}
    \prod_{x \in \Lambda}(1+\tau_x)
    =
    e^{\sum_{x \in \Lambda}\tau_x}
    \prod_{x \in \Lambda}(1+\tau_x)e^{-\tau_x}
    \approx
    e^{\sum_{x \in \Lambda}\tau_x}
    \prod_{x \in \Lambda}e^{-\frac 12 \tau_x^2}.
\end{equation}
The approximation amounts to matching terms up to order $\tau_x^2$ in
a Taylor expansion.
With this approximation, \eqref{e:Gsawid2} corresponds to
\eqref{e:Gsawrep2} with $g=\frac 12$ and $\lambda = -1$.
A careful comparison of the two models is given in
\cite{BS10}.


\section{Supersymmetry}
\label{sec:susy}

%

Integrals such as
$\int e^{-S_A} F(\tau)$
are unchanged if we formally interchange the pairs $\phi,\phib$ and
$\psi,\psib$.  By \eqref{e:F0}, it is also true that
$\int e^{-S_A} F(\tau)\phib_a\phi_b = \int e^{-S_A} F(\tau)\psib_a\psi_b$ (the difference is
$\int e^{-S_A}\tau F(\tau) = 0$).
This suggests the existence of a symmetry between bosons and fermions.
Such a symmetry is called a \emph{supersymmetry}.

In this section, as a brief illustration,
we use methods of supersymmetry to provide
an alternate proof of \eqref{e:F0}, following \cite{BI03d}.
The supersymmetry generator $Q$ is a map on the space of forms which maps
bosons to fermions and vice versa.  It can be defined in
terms of standard operations in differential geometry,
namely the exterior derivative and interior product, as follows.

An antiderivation $F$ is a linear map on forms which obeys
$F(\omega_1\wedge \omega_2)
= F\omega_1\wedge\omega_2 + (-1)^{p_1} \omega_1
\wedge F\omega_2$, when $\omega_1$ is a form of degree $p_1$.
The {\em exterior derivative} $d$ is the linear antiderivation that maps
a form of degree $p$ to a form of degree $p+1$, defined by $d^2=0$ and,
for a zero form $f$,
\begin{equation}
    df = \sum_{x\in \Lambda} \Big(\frac{\partial f}{\partial \phi_x} d\phi_x
    + \frac{\partial f}{\partial \phib_x} d\phib_x \Big).
\end{equation}

Consider the flow
acting on $\C^{\carL}$ defined by
$\phi_x \mapsto e^{ -2 \pi i \theta }\phi_x$.
This flow is generated by the vector field $X$ defined by $X(\phi_x) = -2
\pi i \phi_x$, and $X(\phib_x) = 2 \pi i \phib_x$.  The action by
pullback of the flow on forms is
\begin{equation}
d\phi_x
\mapsto d (e^{-2 \pi i \theta}\phi_x)
 = e^{- 2 \pi i \theta} \, d\phi_x,  \quad \quad
d\phib_x \mapsto e^{ 2 \pi i \theta}\, d\phib_x.
\end{equation}
The {\em interior
product} ${\ci} = \ci_X$ with the vector field $X$ is the
linear antiderivation
that maps forms of degree $p$ to forms of degree $p-1$ (and maps forms
of degree zero to zero), given by
\begin{equation}
{\ci} d\phi_x = -2 \pi i \phi_x,  \quad\quad {\ci} d\phib_x = 2 \pi
i \phib_x.
\end{equation}
The interior product obeys $\ci^2=0$.
\eject

The \emph{supersymmetry generator} $Q$ is defined by
\begin{equation}
 {Q} = d + \ci .
\end{equation}
A form $\omega$ that satisfies $Q \omega = 0$ is called \emph{supersymmetric}
or $Q$-\emph{closed}.   A form $\omega$ that is in the image of $Q$ is called
$Q$-\emph{exact}.
Note that the integral of any $Q$-\emph{exact} form is zero
(assuming that the form decays appropriately
at infinity), since integration acts only on
forms of top degree $2N$ and the degree of
$\ci \omega$ is at most $2N - 1$, while $\int d \omega = 0$ by Stokes'
theorem.
We will use the fact that $Q$ obeys the chain rule for
even forms, in the sense that if
$K = (K_1,\ldots,K_t)$ with each $K_i$ an even form,
and if $F : \C^t \to \C$ is $C^\infty$, then
\begin{equation}
\label{e:Qcr}
    QF(K) = \sum_{i=1}^t F_i(K)
    QK_i,
\end{equation}
where $F_i$ denotes the partial derivative.  A proof is given below.

The \emph{Lie derivative} $\cL = \cL_X$ is the infinitesimal flow
obtained by differentiating with respect to the flow at $\theta = 0$.
Thus, for example,
\begin{equation}
 {\cL} \,d\phi_x = \frac{d}{d\theta} e^{-2\pi i \theta} d\phi_x
  \big|_{\theta =0} = - 2 \pi i \,d\phi_x .
\end{equation}
A form $\omega$ is defined to be \emph{invariant} if
$\cL \omega = 0$.
For example, the form
\begin{equation}
\label{e:uxydef}
    u_{x,y}=\phi_x d\phib_y
\end{equation}
is invariant since it is constant under the flow of $X$.
Cartan's formula asserts that $\cL = d \,\ci + \ci \,d$
(see, e.g., \cite[p.~146]{GHV72}).
Since $d^{2}=0$ and $\ci^2=0$, we have that
$\cL = Q^2$, so $Q$ is the square root of $\cL$.


\smallskip \noindent
\emph{Alternate proof of \eqref{e:F0}.}
We will show that $\int e^{-S_A}F(\lambda \tau)$ is independent of $\lambda \in
\R$.  Comparing the value of this integral for $\lambda =0$ and $\lambda =1$,
the identity \eqref{e:F0} then follows from \eqref{e:sn}.
Computation of the derivative gives
\begin{equation}
\label{e:dlam}
    \frac{d}{d\lambda} \int e^{-S_A}F(\lambda \tau)
    =
    \int e^{-S_A}\sum_{x\in\Lambda}
    F_{x}(\lambda \tau) \tau_x,
\end{equation}
where $F_x$ denotes the partial derivative of $F$ with respect to
coordinate $x$. To show that the integral on the right-hand side
vanishes, it suffices to show that the integrand is $Q$-exact. Let
$v_{x,y}=\frac{1}{2\pi i} u_{x,y}$, where $u_{x,y}$ is given by
\eqref{e:uxydef}. Then $v_{x,y}$ is invariant, and since $Q v_{x,x}
= \tau_{x}$, $\tau_x$ is both $Q$-exact and $Q$-closed. Since
$Q(\sum_{x,y} A_{x,y} v_{x,y}) = S_A$ and $\sum_{x,y} A_{x,y}
v_{x,y}$ is invariant, the form $S_A$ is also $Q$-exact and
$Q$-closed. By \eqref{e:Qcr}, $e^{-S_A}$ and $F_{x}(\lambda \tau)$
are both $Q$-closed. Therefore, since $Q$ is an antiderivation,
\begin{equation}
\label{e:dlam1}
    e^{-S_A}F_x(\lambda \tau)  \tau_x
    =
     Q\left(e^{-S_A} F_{x}(\lambda \tau)
    v_{x,x} \right),
\end{equation}
as required.
\qed

\smallskip \noindent
\emph{Proof of the chain rule \eqref{e:Qcr} for $Q$.}
Suppose first that $K$ is a zero form.  Then
\begin{equation}
    QF(K) = dF(K) = \sum_{i=1}^t
    \left[\frac{\partial F(K)}{\partial \phi_i} d\phi_i +
    \frac{\partial F(K)}{\partial \bar\phi_i} d\bar\phi_i \right].
\end{equation}
By the chain rule, this is $\sum_{i} F_{i}(K) dK_{i}
=\sum_{i} F_{i}(K) QK_{i}$.  This proves \eqref{e:Qcr}
for zero forms, so we may assume now that $K$ is higher degree.

Let $\epsilon_i$ be the multi-index that has $i^{\rm th}$ component $1$ and
all other components $0$.  Let $K^{(0)}$ denote the degree zero part of $K$.
By \eqref{e:Fdef}, the fact that $Q$ is an antiderivation, and the chain
rule applied to zero forms,
\eqalign
    QF(K) & =
    \sum_{\alpha} \frac{1}{\alpha !}
    QF^{(\alpha)}(K^{(0)})
    (K - K^{(0)})^{\alpha}
    +
    \sum_{\alpha} \frac{1}{\alpha !}
    F^{(\alpha)}(K^{(0)})
    Q(K - K^{(0)})^{\alpha}
    \nnb
\label{e:Qcr1}
    & = \sum_{\alpha} \frac{1}{\alpha !}
    \sum_{i=1}^t
    F^{(\alpha + \epsilon_i)}(K^{(0)})[ QK^{(0)}_i]
    (K - K^{(0)})^{\alpha} \nnb
    &\quad\quad+
    \sum_{\alpha} \frac{1}{\alpha !}
    F^{(\alpha)}(K^{(0)})
    Q(K - K^{(0)})^{\alpha}.
\enalign
Since $Q$ is an antiderivation,
\begin{equation}
\label{e:Qcr2}
    Q(K - K^{(0)})^\alpha
    =
    \sum_{i=1}^t
    \alpha_i (K - K^{(0)})^{\alpha -\epsilon_i}[QK_i -
    QK^{(0)}_i].
\end{equation}
The first term on the right-hand side of \eqref{e:Qcr1} is canceled
by the contribution to the second term of \eqref{e:Qcr1} due to the
second term of \eqref{e:Qcr2}. And the contribution to the second
term of \eqref{e:Qcr1} due to the first term of \eqref{e:Qcr2} is
$\sum_{i} F_i(K) QK_{i}$, as required. \qed

\section{Conclusion}

We have given a unified treatment of three representations for
simple random walk in Theorems~\ref{thm:srw}, \ref{thm:srw2} and
\ref{thm:srw3}.  These representations had appeared previously
in \cite{BFS82,Dynk83,BEI92}.  In Theorem~\ref{thm:Gloop}, we
have represented a model of a self-avoiding walk
in a background of self-avoiding loops, all mutually avoiding,
in terms of a (bosonic) Gaussian integral.

Mixed bosonic-fermionic Gaussian integrals were introduced in
Section~\ref{sec:rwdf}, and some elements of the theory of these
integrals were derived. Using these integrals, and particularly
using Proposition~\ref{prop:Ftau0}, representations for the weakly
self-avoiding walk and strictly self-avoiding walk were obtained in
Theorems~\ref{thm:wsaw} and \ref{thm:saw}, respectively. Our
representation in Theorem~\ref{thm:saw} is new. These
representations provide the point of departure for rigorous
renormalization group analyses of various self-avoiding walk
problems \cite{BEI92,BI03c,BI03d,BS10,MS08}. For the strictly
self-avoiding walk, two different proofs of the representation were
given, in Sections~\ref{sec:expres} and \ref{sec:ibp}. The role of
the fermionic part of the representation in eliminating loops was
detailed in Section~\ref{sec:expres}.  This contrasts with the
formal $N \to 0$ limit discussed in Section~\ref{sec:N0lim}.

The mixed bosonic-fermionic representations are examples of
supersymmetric field theories.  A brief discussion of some elements
of supersymmetry was given in Section~\ref{sec:susy}.

\bibliographystyle{acmtrans-ims}

\end{document}